The Annals of Applied Probability
2009, Vol. 19, No. 4, 1385–1403
DOI: 10.1214/08-AAP578
© Institute of Mathematical Statistics, 2009

# ON CONVERGENCE TO STATIONARITY OF FRACTIONAL BROWNIAN STORAGE

By Michel Mandjes, Ilkka Norros and Peter Glynn

*Korteweg-de Vries Institute for Mathematics and University of Amsterdam, VTT Information Technology and Stanford University*

With $M(t) := \sup_{s \in [0,t]} A(s) - s$ denoting the running maximum of a fractional Brownian motion $A(\cdot)$ with negative drift, this paper studies the rate of convergence of $\mathbb{P}(M(t) > x)$ to $\mathbb{P}(M > x)$. We define two metrics that measure the distance between the (complementary) distribution functions $\mathbb{P}(M(t) > \cdot)$ and $\mathbb{P}(M > \cdot)$. Our main result states that both metrics roughly decay as $\exp(-\vartheta t^{2-2H})$, where $\vartheta$ is the decay rate corresponding to the tail distribution of the busy period in an fBm-driven queue, which was computed recently [*Stochastic Process. Appl.* (2006) **116** 1269–1293]. The proofs extensively rely on application of the well-known large deviations theorem for Gaussian processes. We also show that the identified relation between the decay of the convergence metrics and busy-period asymptotics holds in other settings as well, most notably when Gärtner–Ellis-type conditions are fulfilled.

**1. Introduction.** Let $(A(s))_{s \in \mathbb{R}}$ be a centered *fractional Brownian motion* (fBm), that is, a stochastic process, such that for all $s \in \mathbb{R}$, $A(s)$ obeys a Gaussian distribution with mean 0 and variance $|s|^{2H}$, for $H \in (0,1)$. fBm has recently become one of the key models in the applied probability literature, because of a number of interesting features. For $H \in (\frac{1}{2}, 1)$, correlations decay so slowly that the process qualifies as *long-range dependent*; choosing $H \in (\frac{1}{2}, 1)$ leads to positive correlations, whereas $H \in (0, \frac{1}{2})$ results in negative correlations. Also, fBm exhibits *selfsimilar* behavior, in that $A(\alpha s)$ has the same distribution as $\alpha^H A(s)$. Its use has been advocated in several practical settings; see, for example, [21, 28] for networking applications; for applications in finance, see, for example, [2].

Motivated by these applications, substantial attention was paid to the analysis of *regulated fBm*, or *fractional Brownian storage* [27]. With the storage process defined through

$$Q(t) := \sup_{s \le t} A(t) - A(s) - (t - s),$$

distributional properties of the steady-state storage level $Q := \lim_{t \to \infty} Q(t)$ can be used to describe the performance of a network element. The stationary storage level $Q$ being distributed as $M := \sup_{s>0} A(s) - s$, see for instance [22], Sections 5.1 and 5.2, a considerable amount of work [12, 16, 25, 26] has been devoted to the characterization of the distribution of $M$, that is, the supremum attained by an fBm with negative drift. The results obtained are predominantly asymptotic in nature; most notably, an explicit function $\varphi(\cdot)$ was identified [16] such that

$$\mathbb{P}(M > x)/\varphi(x) \to 1$$

as $x \to \infty$. Clearly, $\varphi(x)$ can serve as an approximation of $\mathbb{P}(M > x)$ for large $x$; we lack, however, accurate approximations or bounds for small or moderate values of $x$.

As mentioned above, asymptotic results are available that describe $\mathbb{P}(M > x)$ for $x$ large, but considerably less is known about the convergence of the running supremum

$$M(t) := \sup_{s \in [0,t]} A(s) - s$$

to its limiting distribution $M$. The primary goal of the present paper is to determine the speed of this convergence. Knowledge of this speed of convergence is useful in several contexts. In the first place it provides information on the question whether at the time scale at hand, that is, $t$, it is justified to approximate $\mathbb{P}(M(t) > x)$ by $\mathbb{P}(M > x)$. Loosely speaking, if the convergence to the limiting distribution is fast (i.e., a short relaxation time), then such a procedure gives an accurate approximation, whereas in case of slow convergence transient results need to be used. Supposing that one wishes to estimate $\mathbb{P}(M > x)$ [and hence also $\mathbb{P}(Q > x)$] by performing a stochastic simulation, the speed of convergence also provides us with useful insights with respect to the design of this simulation experiment. It, for instance, sheds light on the question of whether it is more efficient (in terms of the variance of the resulting estimator) to have one long run or multiple shorter runs. In addition, it yields guidelines as to how long one should simulate $A(s) - s$ in order to be able to accurately estimate $\mathbb{P}(M > x)$, that is, to find an appropriate simulation horizon $t$ such that $\mathbb{P}(M(t) > x)$ approximates $\mathbb{P}(M > x)$ sufficiently well (in terms of some explicitly specified criterion); cf. [11].

In principle, the rate of convergence (as $t$ grows large) of $\mathbb{P}(M(t) > x)$ to $\mathbb{P}(M > x)$ could depend on the value of $x$. Obviously, one could use several distance measures, each incorporating this dependence on $x$ in a specific way. First notice that

$$\mathbb{P}(M > x) - \mathbb{P}(M(t) > x) = \mathbb{P}(M > x, M(t) \leq x) =: \gamma(x, t) > 0,$$

as the event $\{M(t) > x\}$ implies $\{M > x\}$. Two possible distances are the following:

- *Kolmogorov–Smirnov distance* (sup-norm). We define (for random variables $X$ and $Y$) $d_1(X, Y)$ by $\sup_x |\mathbb{P}(X \leq x) - \mathbb{P}(Y \leq x)|$. It is well known that $d_1$ is a distance. The first metric used in our paper is based on the distance $d_1$:

$$\mathbb{D}_1(t) := d_1(M, M(t)) = \sup_{x > 0} \gamma(x, t).$$

  $\mathbb{D}_1(t)$ measures the *maximum* distance between the distribution functions of $M$ and $M(t)$.

- *Integral distance* ($L_1$-norm). It is a well-known fact that $d_2(X, Y) := \int_x |\mathbb{P}(X \leq x) - \mathbb{P}(Y \leq x)| \, dx$ is a distance, too. The second metric considered in our paper is based on the distance $d_2$:

$$\mathbb{D}_2(t) := d_2(M, M(t)) = \int_{x > 0} \gamma(x, t) \, dx.$$

  $\mathbb{D}_2(t)$ measures the *total* distance between the distribution functions. Note that $\mathbb{D}_2(t)$ can be interpreted as $\mathbb{E}M - \mathbb{E}M(t)$.

The goal of the paper is to identify the asymptotics of the $\mathbb{D}_i(t)$ for $t$ large, $i = 1, 2$. Our main result is that the decay rates of both distance measures coincide, and are equal to asymptotics of the busy-period distribution in an fBm-driven queue, which were recently identified in [23].

The structure of this paper is as follows. In Section 2 we recapitulate a series of results on the large deviations of fBm, most notably (the generalized version of) Schilder's theorem. We also recall the main results on busy-period asymptotics [23], which enable us to state the main results of our paper. Section 3 presents a number of auxiliary results that are used in Section 4 (in order to determine the asymptotics of the Kolmogorov–Smirnov distance) and Section 5 (in order to determine the asymptotics of the integral distance). In Section 6 we consider the situation of short-range dependent input (or, more precisely, the situation in which so-called Gärtner–Ellis-type conditions are met), to show that also in this regime the asymptotics of the $\mathbb{D}_i(t)$ are equal to those of the busy-period distribution. Section 7 concludes; it includes a procedure for determining the simulation horizon.

**2. Preliminaries and main results.** In this section we recall a number of useful results from the literature. Emphasis is on busy-period asymptotics recently identified in [23]. We then state our main results.

2.1. *Generalized Schilder.* Informally, the generalized version of Schilder's theorem provides us with a "rate functional" $\mathbb{I}(\cdot)$ such that

$$p_n[\mathscr{S}] := \mathbb{P}\left(\frac{A(\cdot)}{\sqrt{n}} \in \mathscr{S}\right) \approx \exp\left(-n \inf_{f \in \mathscr{S}} \mathbb{I}(f)\right).$$

In other words: in this large-deviations setting, the probability of interest decays exponentially in $n$. The "$\approx$" in the above statement should be interpreted as follows: under mild conditions on the set $\mathscr{S}$ (more concretely, if $\mathscr{S}$ is an $\mathbb{I}$-continuity set), the decay rate of $p_n[\mathscr{S}]$ is given by

$$\lim_{n \to \infty} \frac{1}{n} \log p_n[\mathscr{S}] = - \inf_{f \in \mathscr{S}} \mathbb{I}(f).$$

Apart from the Brownian case $H = \frac{1}{2}$, the "rate functional" $\mathbb{I}(\cdot)$ cannot be given explicitly. It is defined through

$$\mathbb{I}(f) := \begin{cases} \frac{1}{2}\|f\|^2, & \text{if } f \in \mathscr{R}; \\ \infty, & \text{otherwise,} \end{cases}$$

where $\mathscr{R}$ is the *reproducing kernel Hilbert space* related to the process $A(\cdot)$— see for details [1, 9]. Here $\|f\| := \sqrt{\langle f, f \rangle}$, where $\langle f, g \rangle$ is a suitably defined inner product between $f, g \in \mathscr{R}$. It is noted that $\mathbb{I}(f)$ can be interpreted as a measure for the "likelihood" of a path $f$; the path $f \equiv 0$ is the sole path that gives $\mathbb{I}(f) = 0$, while for other paths $f$ the rate $\mathbb{I}(f)$ is strictly positive. The path $f^\star := \arg\inf_{f \in \mathscr{S}} \mathbb{I}(f)$ is usually called the *most likely path* in large deviations literature, and it has the interpretation that, conditional on the fBm being in the set $\mathscr{S}$, with overwhelming probability it will be close to $f^\star$; cf. [10].

Later in this paper we repeatedly use the following property. Suppose $f \in \mathscr{R}$, and $g$ is defined by $g(r) = \alpha f(\beta r)$, for $\alpha, \beta > 0$. Then

(1) $$\|g\| = \alpha \beta^H \|f\|.$$

For the purposes of the present paper, more background on "generalized Schilder" is not required; see for a complete account [1, 9, 23].

2.2. *Busy-period asymptotics.* In [23], and its predecessor [29], the focus was on computing the asymptotics, for large values of $t$, of $\mathbb{P}(K > t)$, where

$$K := \inf\{t \geq 0 : Q(t) = 0\} - \sup\{t \leq 0 : Q(t) = 0\}$$

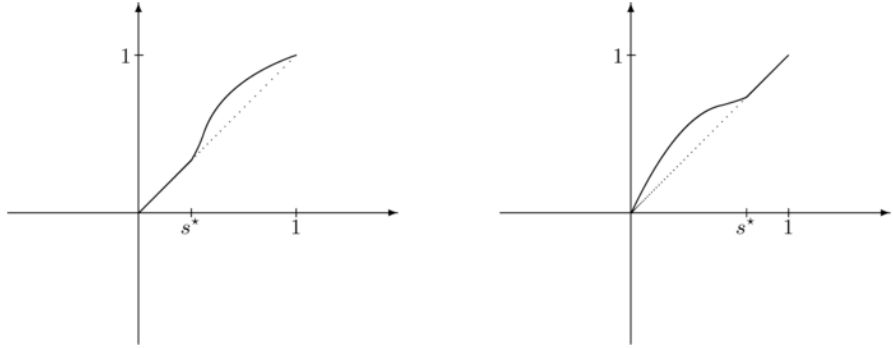

FIG. 1. *Most likely path in $\mathscr{B}$; left picture $H > \frac{1}{2}$; right picture $H < \frac{1}{2}$.*

is the ongoing busy period at time 0. In [29] it was shown that $\mathbb{P}(K > t)$ decays roughly in a Weibullian way, that is, as $\exp(-\vartheta t^{2-2H})$ for some positive constant $\vartheta$. More precisely, it obeys the following logarithmic asymptotics:

$$(2) \qquad \lim_{t \to \infty} \frac{1}{t^{2-2H}} \log \mathbb{P}(K > t) = -\vartheta \qquad \text{where } \vartheta := \inf_{f \in \mathscr{B}} \mathbb{I}(f).$$

Here $\mathscr{B}$ is the set of paths that remain above the diagonal on the interval $[0, 1]$:

$$\mathscr{B} := \{ f \in \mathscr{R} \mid \forall r \in [0, 1] \colon f(r) \geq r \}.$$

Note that the set $\mathscr{B}$ can be regarded as the set of feasible paths that correspond to an *intersection* of events (reflected by the $\forall$-quantor). Where unions are usually easy to deal with, finding the minimizing path in an intersection is typically hard (although often rather precise bounds can be found; see for instance [24]). In [23] we succeeded in determining the right-hand side of (2), as well as the corresponding minimizing path $f^\star \in \mathscr{B}$.

In this particular setting, the most likely paths turn out to have a remarkable shape. For $H > \frac{1}{2}$, the most likely path is at the diagonal in some interval $[0, s^\star]$, and also at time 1, but strictly above the diagonal in between. For $H < \frac{1}{2}$, the corresponding path departs immediately after time 0 from the diagonal, but returns to it strictly before time 1 and continues along it until time 1—see Figure 1. The corresponding decay rate $\vartheta$ is given in [23], Theorem 24.

2.3. *Main results.* We now present the main results of this paper. They entail that both $\mathbb{D}_1(t)$ and $\mathbb{D}_2(t)$ decay as the probability $\mathbb{P}(K > t)$ of the busy period exceeding $t$.

THEOREM 2.1. (i) *For the Kolmogorov–Smirnov distance we have*

$$\lim_{t \to \infty} \frac{1}{t^{2-2H}} \log \mathbb{D}_1(t) = -\inf_{f \in \mathscr{B}} \mathbb{I}(f) = -\vartheta.$$

(ii) *For the integral distance we have*

$$\lim_{t\to\infty} \frac{1}{t^{2-2H}} \log \mathbb{D}_2(t) = -\inf_{f\in\mathscr{B}} \mathbb{I}(f) = -\vartheta.$$

Part (i) and (ii) of this theorem will be proven in Sections 4 and 5, respectively, whereas Section 3 develops a number of useful tools.

**3. Auxiliary results.** In this section we derive a number of results that are needed to prove Theorem 2.1. The following alternative expression for $\gamma(x,t)$ turns out to be useful.

LEMMA 3.1. *For any $x > 0$, $t \geq 0$,*

$$\gamma(x,t) = \mathbb{P}(\forall r \in [0,t] : A(r) \leq x + r; \exists s > t : A(s) > x + s)$$

$$= \mathbb{P}\left(\forall r \in [0,1] : \frac{A(r)}{t^{1-H}} \leq \frac{x}{t} + r; \exists s > 1 : \frac{A(s)}{t^{1-H}} > \frac{x}{t} + s\right).$$

PROOF. The first equality is a matter of rewriting $\{M > x, M(t) \leq x\}$ in terms of the process $A(\cdot)$. The second equality follows from the self-similarity. $\square$

In the sequel we extensively use the following sequence of probabilities:

$$p_n(\delta) := \mathbb{P}\left(\forall r \in [0,1] : \frac{A(r)}{\sqrt{n}} \leq \delta + r; \exists s > 1 : \frac{A(s)}{\sqrt{n}} > \delta + s\right)$$

for $\delta > 0$. We also define their exponential decay rates by

$$J(\delta) := \lim_{n\to\infty} \frac{1}{n} \log p_n(\delta),$$

again for $\delta > 0$; here $J(0)$ denotes the limit of $J(\delta)$ for $\delta \downarrow 0$.

Define $\mathscr{A}_\delta$ as the paths $f$ in the set

(3)     $\mathscr{A}_\delta := \{f \in \mathscr{R} \mid \forall r \in [0,1] : f(r) \leq \delta + r; \exists s > 1 : f(s) > \delta + s\};$

also $\mathscr{A} := \mathscr{A}_0$. The proof of the following result is an immediate consequence of (the generalized version of) Schilder's theorem [1, 9].

LEMMA 3.2. *For any $\delta \geq 0$,*

(4)                              $J(\delta) = -\inf_{f\in\mathscr{A}_\delta} \mathbb{I}(f).$

Define

$$\bar{\mathscr{A}} := \{f \in \mathscr{R} \mid \forall r \in [0,1] : f(r) \leq r; f(1) = 1\};$$

$$\mathscr{B} := \{f \in \mathscr{R} \mid \forall r \in [0,1] : f(r) \geq r\}.$$

The following result concerns a translation of $J(0)$ in terms of our previous result on busy periods, as mentioned in Section 2.2.

Proposition 3.3.

$$\inf_{f \in \mathscr{A}} \mathbb{I}(f) = \inf_{f \in \bar{\mathscr{A}}} \mathbb{I}(f) = \inf_{f \in \mathscr{B}} \mathbb{I}(f).$$

Proof. Due to continuity arguments (cf. the proofs in Section 4 of [29]) the decay rate corresponding to the most likely path in $\mathscr{A}$ is the same as that of the most likely path in

$$\{f \in \mathscr{R} \mid \forall r \in [0,1] : f(r) \leq r; \exists s \geq 1 : f(s) \geq s\}.$$

Let $s^\star$ be the smallest $s \geq 1$ such that $f(s) \geq s$. Define the path $\bar{f}$ through $\bar{f}(r) := f(rs^\star)/s^\star$. Then $\bar{f}(1) = 1$ and, due to (1),

$$\|\bar{f}\|^2 = (s^\star)^{2H-2}\|f\|^2 \leq \|f\|^2.$$

This implies the first equality.

A time shift argument trivially gives that $\inf_{f \in \bar{\mathscr{A}}} \mathbb{I}(f)$ is equal to $\inf_{f \in \bar{\mathscr{A}}_-} \mathbb{I}(f)$, with

$$\bar{\mathscr{A}}_- := \{f \in \mathscr{R} \mid \forall r \in [-1,0] : f(r) \leq r; f(-1) = -1\}.$$

Now reverse time, and we obtain that this infimum is also equal to $\inf_{f \in \bar{\mathscr{B}}} \mathbb{I}(f)$, with

$$\bar{\mathscr{B}} := \{f \in \mathscr{R} \mid \forall r \in [0,1] : f(r) \geq r; f(1) = 1\}.$$

The analysis of [23] implies that the infima over $\mathscr{B}$ and $\bar{\mathscr{B}}$ coincide, which proves the stated. □

To obtain a better intuitive understanding of Proposition 3.2, the reader may compare the most likely busy-period path (i.e., the "cheapest path" in $\mathscr{B}$), as depicted in Figure 1, with the most likely path in $\mathscr{A} = \mathscr{A}_0$, as depicted in Figure 2.

Lemma 3.4. (i) With, for $\delta > 0$,

(5) $$\mathscr{D}_\delta := \{f \in \mathscr{R} \mid \forall r \in [0,1] : f(r) \geq r; f(1) = 1 + \delta\},$$

$\inf_{f \in \mathscr{D}_\delta} \|f\|^2$ increases in $\delta$, for $\delta \in [0, H^{-1} - 1]$.

(ii) With, for $0 < \varepsilon < 1$ and $\delta > 0$,

$$\mathscr{D}_{\delta,\varepsilon} := \{f \in \mathscr{R} \mid \forall r \in [0,1] : f(r) \geq r - \varepsilon; f(1) = 1 + \delta - \varepsilon\},$$

$\inf_{f \in \mathscr{D}_{\delta,\varepsilon}} \|f\|^2$ increases in $\delta$, for $\delta \in [0, H^{-1} - 1]$.

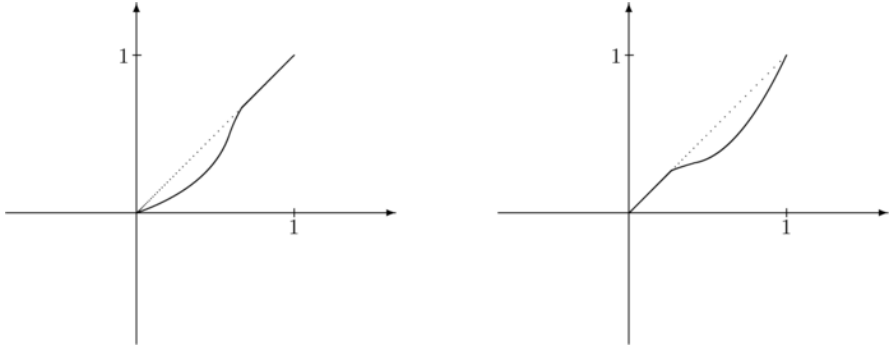

Fig. 2.   *Most likely path in $\mathscr{A}$; left picture $H > \frac{1}{2}$; right picture $H < \frac{1}{2}$.*

Proof.   In order to settle claim (i), first observe that the probability

$$q_n(\delta) := \mathbb{P}\left(\forall r \in [0,1]: \frac{A(r)}{\sqrt{n}} \geq r; \exists s \geq 1: \frac{A(s)}{\sqrt{n}} \geq \delta + s\right),$$

decreases in $\delta$, so consequently also

$$(6) \qquad \lim_{n \to \infty} \frac{1}{n} \log q_n(\delta) = -\frac{1}{2} \inf_{f \in \mathscr{D}_\delta} \|f\|^2$$

decreases in $\delta$; the equality is due to (the generalized version of) Schilder's theorem, and

$$\bar{\mathscr{D}}_\delta := \{f \in \mathscr{R} \mid \forall r \in [0,1]: f(r) \geq r; \exists s \geq 1: f(s) \geq \delta + s\}.$$

Pick an arbitrary path $f$ in this set, and let $s^\star$ be the smallest $s \geq 1$ such that $f(s) \geq s + \delta$. Then consider the path $\bar{f}$, defined by

$$\bar{f}(r) := \left(\frac{1+\delta}{s^\star + \delta}\right) f(rs^\star).$$

Note that $\bar{f}(1) = 1 + \delta$, and, because $f$ lies in the set $\bar{\mathscr{D}}_\delta$, for all $r \in [0,1]$,

$$\bar{f}(r) \geq \left(\frac{1+\delta}{s^\star + \delta}\right) rs^\star;$$

it is easily verified that the right-hand side of the previous display is at least $r$. We conclude that $\bar{f}$ is in $\mathscr{D}_\delta$ as well. Moreover, as before, any path in $\bar{\mathscr{D}}_\delta$ can be replaced by a path in $\mathscr{D}_\delta$ with a smaller norm: due to (1),

$$\|\bar{f}\|^2 = \left(\frac{1+\delta}{s^\star + \delta}\right)^2 (s^\star)^{2H} \|f\|^2 \leq \|f\|^2.$$

Here it is used that for $\delta \in [0, H^{-1} - 1]$ and $s \geq 1$ it holds that $(1+\delta)s^H \leq s + \delta$. In other words, for these $\delta$ we can replace the set $\bar{\mathscr{D}}_\delta$ in (6) by $\mathscr{D}_\delta$.

Hence

$$-\frac{1}{2}\inf_{f\in\mathscr{D}_\delta}\|f\|^2$$

decreases in $\delta$ as well, which implies claim (i).

The proof of claim (ii) is similar. The probability

$$q_n(\delta,\varepsilon):=\mathbb{P}\left(\forall r\in[0,1]:\frac{A(r)}{\sqrt{n}}\geq r-\varepsilon;\exists s\geq 1:\frac{A(s)}{\sqrt{n}}\geq\delta+s-\varepsilon\right),$$

decreases in $\delta$. Again, for any $f$ in

$$\bar{\mathscr{D}}_{\delta,\varepsilon}:=\{f\in\mathscr{R}\mid\forall r\in[0,1]:f(r)\geq r-\varepsilon;\exists s\geq 1:f(s)\geq\delta+s-\varepsilon\},$$

and $s^\star$ as defined before, we can define

$$\bar{f}(r):=\left(\frac{1+\delta-\varepsilon}{s^\star+\delta-\varepsilon}\right)f(rs^\star),$$

which has a smaller norm than $f$ for any $\delta\in[0,H^{-1}-1]$, and which lies in $\mathscr{D}_{\delta,\varepsilon}$; the latter statement follows from $\bar{f}(1)=1+\delta-\varepsilon$, in combination with, for $r\in[0,1]$,

$$\bar{f}(r)\geq\left(\frac{1+\delta-\varepsilon}{s^\star+\delta-\varepsilon}\right)(rs^\star-\varepsilon)\geq r-\varepsilon,$$

where the first inequality is due to $f\in\bar{\mathscr{D}}_{\delta,\varepsilon}$, and the second inequality due to the fact that, obviously, $(s^\star-1)(\delta r+\varepsilon(1-r))\geq 0$ (realize that $1+\delta-\varepsilon>0$ because $\varepsilon<1$). We have now established claim (ii). $\square$

We now establish a useful lemma on the behavior of $J(\delta)$, that is, the decay rate of $p_n(\delta)$, as a function of $\delta$.

Lemma 3.5. (i) $J(\delta)$ decreases for $\delta\in(0,H^{-1}-1]$.
(ii) $J(\delta)<J(0)$ for $\delta>H^{-1}-1$.

Proof. First consider part (i). Due to continuity arguments, and "generalized Schilder," the decay rate $J(\delta)$ equals the decay rate corresponding to the most likely path in

(7) $$\{f\in\mathscr{R}\mid\forall r\in[0,1]:f(r)\leq r+\delta;\exists s\geq 1:f(s)\geq s+\delta\}.$$

With arguments identical to those used in the proof of Lemma 3.4, we can show that for $\delta\in[0,H^{-1}-1]$ we can replace the set (7) by

(8) $$\{f\in\mathscr{R}\mid\forall r\in[0,1]:f(r)\leq r+\delta;f(1)=1+\delta\}.$$

Reversing time, we observe that the decay rate corresponding with the most likely path in (8) equals the one corresponding to the set $\mathscr{D}_\delta$, as defined through (5). Observe that the event corresponding to the latter set becomes increasingly rare when $\delta$ grows; Lemma 3.4(i) now implies claim (i).

Now consider the second claim. Clearly

$$J(\delta) \leq \lim_{n \to \infty} \frac{1}{n} \log \mathbb{P}\left(\exists s > 1 : \frac{A(s)}{\sqrt{n}} > \delta + s\right).$$

The decay rate in the right-hand side equals (see, e.g., [22], Exercise 6.1.3)

$$-\inf_{s \geq 1} \frac{(s + \delta)^2}{2s^{2H}},$$

which reduces, for $\delta > H^{-1} - 1$, to

$$(9) \qquad\qquad -\frac{1}{2}\left(\frac{\delta}{1-H}\right)^{2-2H}\left(\frac{1}{H}\right)^{2H};$$

this is a decreasing function, with value $-\frac{1}{2}H^{-2}$ for $\delta = H^{-1} - 1$. Conclude that, for $\delta > H^{-1} - 1$, $J(\delta) \leq -\frac{1}{2}H^{-2}$. It is left to prove that $J(0) > -\frac{1}{2}H^{-2}$. To this end, we use that in [29] it was shown that

$$-\frac{1}{2} \cdot \varphi(H) \leq J(0) \leq -\frac{1}{2};$$

$$\varphi(H) := \frac{1}{H(2H-1)(2-2H)} \cdot \frac{\Gamma(3/2 - H)}{\Gamma(H - 1/2)\Gamma(2 - 2H)}.$$

As seen from Figure 4.1 in [29], we need to check $\varphi(H) < H^{-2}$ for $H$ in the neighborhood of 1. Observe that both functions have value 1 for $H = 1$. As $\varphi(\cdot)$ is concave in this neighborhood, and $H^{-2}$ convex, we need to verify whether $\varphi'(1) > -2$. Calculations yield that, with $\gamma_{\mathrm{EM}}$ the Euler–Mascheroni constant 0.5772,

$$\varphi'(1) = -3 - 2 \cdot \frac{\Gamma'(1/2)}{\Gamma(1/2)} - 2\gamma_{\mathrm{EM}} \approx -0.22 > -2.$$

This proves claim (ii).  □

Lemma 3.5 is already a first indication that Theorem 2.1(i) indeed holds, as seen as follows. Observe that (provided that the limits exist), using Lemma 3.1, with $\alpha$ defined as $(2 - 2H)^{-1}$,

$$\lim_{t \to \infty} \frac{1}{t^{2-2H}} \log \mathbb{D}_1(t)$$

$$= \lim_{t \to \infty} \frac{1}{t^{2-2H}} \log \sup_{x>0} \gamma(x, t)$$

$$= \lim_{n \to \infty} \frac{1}{n} \log \sup_{x>0} \mathbb{P}\left(\forall r \in [0, 1] : \frac{A(r)}{\sqrt{n}} \leq \frac{x}{n^\alpha} + r;\right.$$

$$(10) \qquad\qquad \exists s > 1 : \frac{A(s)}{\sqrt{n}} > \frac{x}{n^{\alpha}} + s \Big)$$

$$= \lim_{n \to \infty} \frac{1}{n} \log \sup_{x > 0} \mathbb{P}\Big( \forall r \in [0,1] : \frac{A(r)}{\sqrt{n}} \le x + r;$$

$$\exists s > 1 : \frac{A(s)}{\sqrt{n}} > x + s \Big)$$

$$= \lim_{n \to \infty} \frac{1}{n} \log \sup_{x > 0} p_n(x).$$

In other words, if one can interchange the limit and supremum, then

$$\lim_{t \to \infty} \frac{1}{t^{2-2H}} \log \mathbb{D}_1(t) = \sup_{x > 0} J(x),$$

which equals $J(0)$ due to Lemma 3.5; applying Lemma 3.2 and Proposition 3.3, conclude that this would also mean that Theorem 2.1(i) holds. The goal of Section 4 is to prove that the limit and supremum can indeed be interchanged.

In light of the fact that, in a large-deviations setting, the decay rate of an integral is, under rather general conditions, determined by the decay rate of the maximum of the integrand, it is now also expected that indeed Theorem 2.1(ii) holds. Section 5 is devoted to substantiating this claim.

## 4. Proof for Kolmogorov–Smirnov distance.

PROOF OF THEOREM 2.1(i). As seen above in (10),

$$\lim_{t \to \infty} \frac{1}{t^{2-2H}} \log \mathbb{D}_1(t) = \lim_{n \to \infty} \frac{1}{n} \log \sup_{x > 0} p_n(x),$$

implying that the lower bound is trivial, as for all $\varepsilon > 0$,

$$\lim_{t \to \infty} \frac{1}{t^{2-2H}} \log \mathbb{D}_1(t) \ge J(\varepsilon).$$

Now let $\varepsilon \downarrow 0$, and apply Lemma 3.2 and Proposition 3.3. The upper bound is proven in the following steps:

(i) First observe that, taking for convenience $M$ as a multiple of $\varepsilon > 0$,

$$\sup_{x > 0} p_n(x) \le \left( \sum_{k=1}^{M/\varepsilon} p_n^{\varepsilon}(k) \right) + \sup_{x > M} p_n(x)$$

$$(11)$$

$$\text{with } p_n^{\varepsilon}(k) := \sup_{x \in ((k-1)\varepsilon, k\varepsilon]} p_n(x).$$

(ii) Trivially,

$$\sup_{x>M} p_n(x) \leq \sup_{x>M} \mathbb{P}\left(\exists s > 1 : \frac{A(s)}{\sqrt{n}} > x + s\right) = \mathbb{P}\left(\exists s > 1 : \frac{A(s)}{\sqrt{n}} > M + s\right).$$

Also, as seen in (9), for $M > H^{-1} - 1$,

$$\lim_{n\to\infty} \frac{1}{n} \log \mathbb{P}\left(\exists s > 1 : \frac{A(s)}{\sqrt{n}} > M + s\right) = -\frac{1}{2}\left(\frac{M}{1-H}\right)^{2-2H}\left(\frac{1}{H}\right)^{2H}.$$

We saw in the proof of Lemma 3.5 that, for $M \geq H^{-1} - 1$, this expression is majorized by $J(0)$, and hence

(12)                    $$\lim_{n\to\infty} \frac{1}{n} \log \sup_{x>M} p_n(x) \leq J(0).$$

From now on we pick $M := H^{-1} - 1$.

(iii) Also,

$$p_n^\varepsilon(k) \leq \mathbb{P}\left(\forall r \in [0,1] : \frac{A(r)}{\sqrt{n}} \leq k\varepsilon + r; \exists s > 1 : \frac{A(s)}{\sqrt{n}} > (k-1)\varepsilon + s\right);$$

because of "generalized Schilder," we thus obtain

$$\lim_{n\to\infty} \frac{1}{n} \log p_n^\varepsilon(k) \leq -\inf_{f \in \mathscr{A}_{k\varepsilon,\varepsilon}} \mathbb{I}(f),$$

where

$$\mathscr{A}_{x,\varepsilon} := \{f \in \mathscr{R} \mid \forall r \in [0,1] : f(r) \leq x + r; \exists s > 1 : f(s) > x - \varepsilon + s\}.$$

Applying Lemma A.1 in the Appendix to the right-hand side of (11), and using (12), we obtain

(13)        $$\lim_{n\to\infty} \frac{1}{n} \log \sup_{x>0} p_n(x) \leq \max\left\{\max_{x=\varepsilon,2\varepsilon,\ldots,M}\left(-\inf_{f\in\mathscr{A}_{x,\varepsilon}} \mathbb{I}(f)\right), J(0)\right\}.$$

(iv) We now show that, for $x = \varepsilon, 2\varepsilon, \ldots, M$,

(14)                    $$-\inf_{f\in\mathscr{A}_{x,\varepsilon}} \mathbb{I}(f) \leq -\inf_{f\in\mathscr{A}_{0,\varepsilon}} \mathbb{I}(f).$$

This is done as in Lemma 3.5(i). First, using the arguments of the proof of Lemma 3.4, we can restrict ourselves for $x \in \{\varepsilon, 2\varepsilon, \ldots, M\}$ to the paths that attain the value $1 + x - \varepsilon$ at time 1, that is, paths in

$$\{f \in \mathscr{R} \mid \forall r \in [0,1] : f(r) \leq x + r; f(1) = 1 + x - \varepsilon\}.$$

Reversing time yields that this is equivalent to finding the most likely path in

$$\{f \in \mathscr{R} \mid \forall r \in [0,1] : f(r) \geq r - \varepsilon; f(1) = 1 + x - \varepsilon\}.$$

We now observe that the event corresponding to this set is increasingly rare for growing $x$, as made precise by Lemma 3.4(ii) (where we have chosen $\varepsilon < 1$). Hence, for $x = \varepsilon, 2\varepsilon, \ldots, M$, claim (14) now follows.

(v) Hence the right-hand side of (13) is bounded from above by

$$(15) \qquad \max\Big\{-\inf_{f\in\mathscr{A}_{0,\varepsilon}}\mathbb{I}(f),J(0)\Big\}.$$

By letting $\varepsilon\downarrow 0$,

$$J_\varepsilon(0):=\Big(-\inf_{f\in\mathscr{A}_{0,\varepsilon}}\mathbb{I}(f)\Big)\downarrow\Big(-\inf_{f\in\mathscr{A}_0}\mathbb{I}(f)\Big)=J(0).$$

Now the stated follows from letting $\varepsilon\downarrow 0$ in (15), and application of Lemma 3.2 and Proposition 3.3. □

## 5. Proof for integral distance.

PROOF OF THEOREM 2.1(ii). We start by establishing the lower bound. Evidently, for $\varepsilon>0$ arbitrarily chosen, and, as before, $\alpha:=(2-2H)^{-1}$,

$$\int_{x>0}\gamma(x,t)\,dx\geq\int_{x\in[\varepsilon n^\alpha,2\varepsilon n^\alpha]}\gamma(x,t)\,dx\geq\varepsilon n^\alpha\Big(\inf_{x\in[\varepsilon n^\alpha,2\varepsilon n^\alpha]}\gamma(x,t)\Big).$$

Hence, for $\varepsilon>0$ arbitrarily small, using $n^{-1}\cdot\log n\to 0$,

$$\liminf_{t\to\infty}\frac{1}{t^{2-2H}}\log\mathbb{D}_2(t)\geq\liminf_{t\to\infty}\frac{1}{n}\log\Big(\inf_{x\in[\varepsilon n^\alpha,2\varepsilon n^\alpha]}\gamma(x,n^\alpha)\Big)$$

$$=\liminf_{n\to\infty}\frac{1}{n}\log\Big(\inf_{x\in[\varepsilon,2\varepsilon]}p_n(x)\Big).$$

It is straightforward that

$$\inf_{x\in[\varepsilon,2\varepsilon]}p_n(x)\geq\mathbb{P}\Big(\forall r\in[0,1]:\frac{A(r)}{\sqrt{n}}\leq\varepsilon+r;\exists s>1:\frac{A(s)}{\sqrt{n}}>2\varepsilon+s\Big).$$

The lower bound now follows from letting $\varepsilon\downarrow 0$, together with the usual continuity arguments and time reversal.

We now turn to the upper bound. Obviously, for any $\varepsilon>0$,

$$\int_{x>0}\gamma(x,n^\alpha)\,dx\leq n^{\alpha+\varepsilon}\Big(\sup_{x\in(0,n^{\alpha+\varepsilon}]}\gamma(x,n^\alpha)\Big)+\sum_{k=\lceil n^{\alpha+\varepsilon}\rceil}^\infty\Big(\sup_{k\in[k,k+1)}\gamma(x,n^\alpha)\Big).$$

We consider the decay rates (in $n$) of both terms. First focus on the first term; because of, again, $n^{-1}\cdot\log n\to 0$, and Theorem 2.1(i),

$$(16) \qquad \limsup_{n\to\infty}\frac{1}{n}\log\Big(n^{\alpha+\varepsilon}\Big(\sup_{x\in(0,n^{\alpha+\varepsilon}]}\gamma(x,n^\alpha)\Big)\Big)$$

$$\leq\limsup_{n\to\infty}\frac{1}{n}\log\Big(\sup_{x\in(0,\infty)}\gamma(x,n^\alpha)\Big)=-\inf_{f\in\mathscr{B}}\mathbb{I}(f).$$

Also, applying part (i) of Lemma A.2 in the Appendix, there are positive constants $\kappa, \lambda$ such that

$$
\begin{aligned}
\limsup_{n\to\infty} \frac{1}{n} &\log\left(\sum_{k=\lfloor n^{\alpha+\varepsilon}\rfloor}^{\infty} \left(\sup_{x\in[k,k+1)} \gamma(x,n^{\alpha})\right)\right) \\
&\leq \limsup_{n\to\infty} \frac{1}{n} \log\left(\sum_{k=\lfloor n^{\alpha+\varepsilon}\rfloor}^{\infty} \left(\sup_{x\in[k,k+1)} \mathbb{P}(M>x)\right)\right) \\
&\leq \limsup_{n\to\infty} \frac{1}{n} \log\left(\sum_{k=\lfloor n^{\alpha+\varepsilon}\rfloor}^{\infty} \left(\sup_{x\in[k,k+1)} \kappa\exp(-\lambda x^{2-2H})\right)\right) \\
&\leq \limsup_{n\to\infty} \frac{1}{n} \log\left(\sum_{k=\lfloor n^{\alpha+\varepsilon}\rfloor}^{\infty} \kappa\exp(-\lambda k^{2-2H})\right).
\end{aligned}
\tag{17}
$$

Now part (ii) of Lemma A.2 can be applied: there exist positive constants $\bar{\kappa}, \bar{\lambda}$ such that (17) is bounded above by

$$
\limsup_{n\to\infty} \frac{1}{n} \log(\bar{\kappa}\exp(-\bar{\lambda}(n^{\alpha+\varepsilon})^{2-2H})) = -\infty,
$$

where the last equality follows by recalling that $\alpha = (2-2H)^{-1}$. Combining this with (16) and using Lemma A.1 completes the upper bound. $\square$

One might think that asymptotic results for $\mathbb{D}_1$ are stronger than those for $\mathbb{D}_2$, as they involve the whole distribution rather than just the first moment. However, the metrics are different and a priori the asymptotics could differ. Consider for instance a situation in which $\mathbb{P}(Z>x) = x^{-\alpha}$ (where $x>1$, $\alpha>1$) and $\mathbb{P}(Z(t)>x) = 1_{(1,t)}x^{-\alpha}$ (again $x>1$, $\alpha>1$), where we have that $d_1(Z, Z(t)) = t^{-\alpha}$ but $d_2(Z, Z(t)) = t^{-\alpha+1}/(\alpha-1)$. It is noted however that, despite the fact that the asymptotics may differ, both notions are strongly related under rather broad circumstances; see, for instance, the proof of Theorem 2.1(ii), where we could make use of the result for the asymptotics of $\mathbb{D}_1(t)$ to establish the corresponding result for the asymptotics of $\mathbb{D}_2(t)$.

**6. Analogous results for short-range dependent input.** The analysis for fBm shows that the logarithmic asymptotics of both distances coincide with those of long busy periods. One may wonder whether such a property is valid under more general circumstances. One could pursue to extend the class of models for which this result holds to Gaussian processes with regularly varying variance functions; cf. [8, 10]. In this section we focus on non-Gaussian processes, namely, short-range dependence processes that obey *Gärtner–Ellis-type conditions*; see, for example, [17].

To this end, with $A(t)$, as before, the traffic generated by a process with stationary increments in a window of length $t$ (which we can assume to have zero mean, without loss of generality), consider for $x > 0$

$$K(x) := \lim_{t \to \infty} \frac{1}{t} \log \mathbb{P}(\forall r \in [0, t] : A(r) \le xt + r; \exists s > t : A(s) \ge xt + s).$$

It is trivial to rewrite this decay rate to

$$\lim_{t \to \infty} \frac{1}{t} \log \mathbb{P}(\forall r \in [0, 1] : t^{-1} \cdot A(rt) \le x + r; \exists s > 1 : t^{-1} \cdot A(st) \ge x + s).$$

First define the (*asymptotic*) *cumulant function*

$$\Lambda(s) := \lim_{t \to \infty} \frac{1}{t} \log \mathbb{E} e^{sA(t)},$$

which we assume to exist; this essentially means that the input traffic is short-range dependent. Then one can define the *large deviations rate function* by its Legendre transform

$$I(a) := \sup_s (sa - \Lambda(s)).$$

It is readily verified that, under mild conditions, the decay rate in the previous display is bounded from above by

$$(18) \qquad \lim_{t \to \infty} \frac{1}{t} \log \mathbb{P}(\exists s > 1 : t^{-1} \cdot A(st) \ge x + s) = -\inf_{s \ge 1} sI\left(\frac{x + s}{s}\right);$$

these mild conditions in particular relate to the behavior of the input process between grid points, as formalized in [12], Hypothesis 2.3.

We now study under which conditions this upper bound is actually tight. First observe that for any $T > 1$,

$$K(x) \ge \lim_{t \to \infty} \frac{1}{t} \log \mathbb{P}(\forall r \in [0, 1] : t^{-1} \cdot A(rt) \le x + r;$$

$$\exists s \in (1, T] : t^{-1} \cdot A(st) \ge x + s);$$

below we specify how $T$ should be chosen. Let the infimum in the right-hand side of (18) be attained in $s^\star \ge 1$. Suppose $t^{-1} \cdot A(\cdot t)$ satisfies a sample-path large deviations principle (sp-LDP) on $[0, T]$ of the Mogulskii type, with rate function $\bar{\bar{I}}(\cdot)$ and supremum-norm, then

$$K(x) \ge -\bar{\bar{I}}(f)$$

for any feasible path $f$, that is, all $f$ in $\mathscr{A}_x$, as defined by (3); here $\bar{\bar{I}}(f) := \int_0^T I(f'(t)) \, dt$. Now verify that the path given by $f^\star(s) = s \cdot (x + s^\star)/s^\star$ for $s \in [0, s^\star]$, and $f(s) = x + s^\star$ for $s > s^\star$, is indeed feasible; also

$$\bar{\bar{I}}(f^\star) = s^\star I\left(\frac{x + s^\star}{s^\star}\right),$$

so that we can conclude that indeed

$$K(x) = -\inf_{s \geq 1} sI\left(\frac{x+s}{s}\right).$$

Now that we have an expression for $K(x)$, we wonder whether we can prove the decay of the $\mathbb{D}_i(t)$ similar to that of the tail distribution of the busy period. To this end, first note that $K(x)$ is a decreasing function of $x$. It is relatively straightforward to check that this fact entails that, for $i = 1, 2$,

$$\lim_{t \to \infty} \frac{1}{t} \log \mathbb{D}_i(t) = K(0) = -I(1).$$

This argumentation indicates that the decay rate of $\mathbb{D}_i(t)$ indeed coincides with the busy-period asymptotics $-I(1)$, like in the fBm case, as long as an sp-LDP is available. For a broad class of discrete-time processes satisfying a Gärtner–Ellis condition (covering discrete-time Markov modulated processes), such an sp-LDP was proved by Chang [6], whereas for Lévy processes see, for example, [7]. It is clear that, to make the above argumentation work, it is sufficient that $T$ is chosen larger than $s^\star$.

## 7. Discussion and concluding remarks.

*Decay of the correlation function.*    In Section 6 we showed for short-range dependent input that, under specific regularity conditions, $\mathbb{D}_i(t)$ $(i = 1, 2)$ decay essentially exponentially in $t$, and this decay roughly coincides with that of the tail of the busy-period distribution. Let us now consider the asymptotics of the covariance $\mathbb{C}\mathrm{ov}(Q(0), Q(t))$, with, as before,

$$Q(t) := \sup_{s \leq t} A(t) - A(s) - (t - s),$$

assuming that the queue is in equilibrium at time 0. In [14] it was shown for short-range dependent Lévy input (with no negative jumps), that also $\mathbb{C}\mathrm{ov}(Q(0), Q(t))$ has the same asymptotic behavior as the tail of the busy-period distribution: for $i = 1, 2$,

$$\lim_{t \to \infty} \frac{1}{t} \log \mathbb{C}\mathrm{ov}(Q(0), Q(t)) = \lim_{t \to \infty} \frac{1}{t} \log \mathbb{D}_i(t) = -I(1)$$

with $I(1)$ as defined in Section 6.

Our main result states that also in the fBm case, we saw that the busy-period asymptotics and those of $\mathbb{D}_i(t)$ $(i = 1, 2)$ match. In light of the above findings for short-range dependent Lévy input, this suggests that also in the fBm case

$$\lim_{t \to \infty} \frac{1}{t^{2-2H}} \log \mathbb{C}\mathrm{ov}(Q(0), Q(t)) = -\vartheta.$$

Based on the recent results in [13], however, we expect that this is *not* true. Instead, we anticipate that the asymptotics of $\mathbb{C}\mathrm{ov}(Q(0), Q(t))$ are roughly polynomially, or, more precisely, decaying as $t^{2H-2}$, which is equally fast as the asymptotics of $\mathbb{C}\mathrm{ov}(A(0, 1), A(t, t+1))$. A formal proof of this property is still lacking, though. As we feel that the determination of the correlation asymptotics of reflected fBm is an important open problem, we state it as a conjecture.

CONJECTURE 7.1. *For some constant $\gamma \in (0, \infty)$,*
$$\lim_{t \to \infty} t^{2-2H} \, \mathbb{C}\mathrm{ov}(Q(0), Q(t)) = \gamma.$$

*Relation with busy-period asymptotics.* The relation between the asymptotics of the metrics considered and those of the tail of the busy period has been observed before in various contexts before. It is noted that in a number of specific Markovian settings one has shown that the rate of convergence does not depend on the initial state, which is a result that is reminiscent of ours; we refer to [4, 5, 30], which are closely related to generic results by Kingman [19, 20], for the case of birth–death processes. Regarding the fact that the decay rate of the $\mathbb{D}_i(t)$ coincides with that of the tail of the busy-period distribution, we mention a classical result by Kingman [18], Theorem 7, Lemma 7, in the setting of an M/G/1 queue, and results for various queueing systems by Blanc and van Doorn [3]. For the special case of spectrally positive Lévy input the double Laplace transform of $\gamma(x, t)$ can be determined in a similar way as in [14], and inversion techniques can be applied to identify the *exact* asymptotics [i.e., a function $\psi_x(\cdot)$ is found such that $\psi_x(t) \cdot \gamma(x, t) \to 1$ as $t \to \infty$]; elementary calculations show that then also the busy-period decay rate appears.

On an intuitive level it is clear that there is a relation between the rate of convergence to stationarity and the tail of the busy-period distribution. In case of short-range dependent input (and in particular when busy periods are independent) the autocorrelation essentially breaks when the busy period ends. Our result shows that this heuristic carries over to the setting of long-range dependent input (in which there is still substantial dependence between consecutive busy periods). We remark however, that under general conditions (i.e., irrespective of the traffic being long-range dependent or short-range dependent) upper bounds on the metric $\mathbb{D}_1(t)$ in terms of the tail distribution of an ongoing busy period, can be constructed relying on a coupling argument [15].

*Use of convergence estimates in simulation.* If one aims at estimating the probability $\mathbb{P}(M > x)$ through simulation, one needs to truncate the infinite-time horizon to some finite value, say $t$. This evidently always implies an underestimation. Obviously, one needs to choose $t$ sufficiently large such that the error made is negligible. Let $T_x$ be defined as the smallest $s$ such that $A(s) - s = x$; then it holds that $\mathbb{P}(M > x) = \mathbb{P}(T_x < \infty)$. Then one could, for instance, require for small $\varepsilon > 0$ (for instance 5%) that $t$ be chosen large enough that

$$\frac{\mathbb{P}(t < T_x < \infty)}{\mathbb{P}(T_x < \infty)} < \varepsilon.$$

The numerator equals $\gamma(x,t)$, and can therefore be approximated by $\exp(-\vartheta \times t^{2-2H})$ (noting that Lemma 3.5 suggests that this is a conservative estimate), whereas the denominator can be bounded from below by [27] $1 - \Phi(x + t^\star)$, with $\Phi(\cdot)$ denoting the standard normal distribution function, and $t^\star := xH/(1 - H)$; $1 - \Phi(x + t^\star)$ is approximated by

$$\exp\left(-\frac{1}{2}\frac{(x + t^\star)^2}{(t^\star)^{2H}}\right) = \exp\left(-\frac{1}{2}\left(\frac{x}{1 - H}\right)^{2-2H}\left(\frac{1}{H}\right)^{2H}\right).$$

In this way we can compute an estimate for the simulation horizon $t$:

$$t \geq \left(-\frac{\log \varepsilon}{\vartheta} + \frac{1}{2\vartheta}\left(\frac{x}{1 - H}\right)^{2-2H}\left(\frac{1}{H}\right)^{2H}\right)^{1/(2-2H)}.$$

The horizon grows in $x$, as expected, and it does so in a linear fashion for $x$ large. A procedure for the discrete-time counterpart is detailed in [11], see also Proposition 8.1.1 in [22]; it shows the same qualitative behavior as a function of $x$. A procedure for numerically computing $\vartheta$ can be found in [23].

As mentioned in the Introduction, a second evident application of our estimates relates to another issue in the design of the above simulation experiment: it provides insight into the question whether it is more efficient to simulate one long run, or to simulate multiple shorter runs.

## APPENDIX: USEFUL BOUNDS

LEMMA A.1.   *Let, for $i$ in some finite index set $I$, $a_n^{(i)}$ be sequences such that*

$$\limsup_{n \to \infty} \frac{1}{n} \log a_n^{(i)} \leq \omega_i. \tag{19}$$

*Then*

$$\limsup_{n \to \infty} \frac{1}{n} \log\left(\sum_{i \in I} a_n^{(i)}\right) \leq \omega^\star := \max_{i \in I} \omega_i.$$

PROOF. Although we believe the proof is rather standard, we present it here. Choose an arbitrary $\varepsilon > 0$. Then (19) entails that there is an $n_i$ such that for all $n > n_i$ we have that $a_n^{(i)} \leq \exp(n(\omega_i + \varepsilon))$. Hence, for $n > \max_i n_i$,

$$a_n^{(i)} \leq \exp(n(\omega^\star + \varepsilon)).$$

Then

$$\limsup_{n \to \infty} \frac{1}{n} \log \left( \sum_{i \in I} a_n^{(i)} \right) \leq \limsup_{n \to \infty} \frac{1}{n} \log(\{\#I\} \times e^{n(\omega^\star + \varepsilon)}) = \omega^\star + \varepsilon.$$

The stated follows after sending $\varepsilon \downarrow 0$. $\quad\square$

LEMMA A.2. (i) *There exist positive constants $\kappa$ and $\lambda$ such that*

$$\mathbb{P}(M > x) \leq \kappa \exp(-\lambda x^{2-2H}).$$

(ii) *There exists positive constants $\bar{\kappa}$ and $\bar{\lambda}$ such that*

$$\sum_{k=K}^{\infty} \exp(-k^{2-2H}) \leq \bar{\kappa} \exp(-\bar{\lambda} K^{2-2H}).$$

PROOF. Part (i) follows immediately from Duffield and O'Connell [12], Section 3.2; part (ii) is due to Dieker and Mandjes [11], Lemmas 2.1 and 2.2. $\quad\square$

CONVERGENCE TO STATIONARITY OF FRACTIONAL BROWNIAN STORAGE 21

M. Mandjes
Korteweg-de Vries Institute
  for Mathematics
Plantage Muidergracht 24
1018 TV Amsterdam
The Netherlands
and
Stanford University
Stanford, California 94305
USA
E-mail: m.r.h.mandjes@uva.nl

I. Norros
VTT Information Technology
P.O. Box 1202, VTT
Finland
E-mail: ilkka.norros@vtt.fi

P. Glynn
Stanford University
Stanford, California 94305
USA
E-mail: glynn@stanford.edu